\documentclass[14pt]{extarticle}
\usepackage[cp1251]{inputenc}
\usepackage{amsmath}
\usepackage{amssymb}
\usepackage[russian,english]{babel} % We're not running pdftex

\usepackage[pdftex]{graphicx}

 \numberwithin{equation}{section}

\renewcommand{\d}{\delta}
\newcommand{\R}{\mathbb{R}}

\newcommand{\Z}{\mathbb{Z}}
\newcommand{\F}{\mathbb{F}}
\newcommand{\e}{\varepsilon}
\renewcommand{\geq}{\geqslant}
\renewcommand{\leq}{\leqslant}
\renewcommand{\l}{\lambda}
\renewcommand{\o}{\omega}
\renewcommand{\u}{\mathbf u}
\renewcommand{\v}{\mathbf v}

\newtheorem{theorem*}{Theorem}%[section]
\newtheorem{lem}{Lemma}

\begin{document}

\begin{center}
	\textbf{ESTIMATES FOR CHARACTER SUMS \\
		IN FINITE FIELDS OF ORDER $p^2$ and $p^3$}
\end{center}

\begin{center}
	\textbf{Mikhail Gabdullin}  \footnote{The work is supported by the grant from Russian Science Foundation (Project 14-11-00702).} \footnote{E-mail: gabdullin.mikhail@yandex.ru}

\smallskip

\textit{Lomonosov Moscow State University \\
	Krasovskii Institute of Mathematics and Mechanics, Yekaterinburg}

\end{center}

%\selectlanguage{english}

\begin{abstract}
	We obtain nontrivial bounds on character sums over ``boxes'' of volume $p^{n(1/4+\e)}$ in finite fields of order $p^n$ for the cases $n=2$ and $n=3$.
\end{abstract}

\section{Introduction}

Let $p$ be a prime number, $\F_{p^n}$ be the finite field of order $p^n$, and $\{\o_1,\ldots\o_n\}$ be a basis of $\F_{p^n}$ over $\F_p$. Let, further, $N_i,H_i$ be integers such that $1\leq H_i\leq p$,  
$\,\,i=1,\ldots,n$. Define $n$-dimensional parallelepiped $B\subseteq\F_{p^n}$ as follows: 
$$B=\left\{\sum_{i=1}^nx_i\o_i \,:\, N_i+1\leq x_i\leq N_i+H_i, \,\,\, 1\leq i\leq n\right\}. 
$$

%\selectlanguage{russian}

We are interested in estimates for sums $\sum_{x\in B}\chi(x)$, where $\chi$ is a nontrivial multiplicative character of $\F_{p^n}$, with the possible weakest restrictions on $B$. First we give a survey of known results in this direction. In the case $n=1$, more than half a century Burgess's estimate \cite{Burg1} remains to be the strongest one: for every $\e>0$ there exists a $\d>0$ such that for all $H\geq p^{1/4+\e}$ the following inequality holds:
\begin{equation*} 
\left|\sum_{x=N+1}^{N+H}\chi(x)\right|\ll_{\e}p^{-\d}H. 
\end{equation*}
Also Burgess \cite{Burg2} proved an analog of this inequality for $n=2$ and special bases and Karatsuba \cite{Kar1}, \cite{Kar2} generalized it for arbitrary finite fields; so, for instance, in \cite{Kar2} the case of basis $\o_i=g^i$ is considered, where $g$ is a root of an irreducible polynom of degree $n$ over $\F_p$. With this connection it looks natural to find estimates which hold uniformly over all bases of $\F_{p^n}$. Davenport and Lewis were the first to obtain such a result \cite{DL}.

\bigskip

\textbf{Theorem A} (\cite{DL}). \textit{For every $\e>0$ there exists a $\d>0$ such that if
$$H_1=\ldots=H_n=H>p^{\frac{n}{2n+2}+\e},
$$
then}
$$\left|\sum_{x\in B}\chi(x)\right|\leq(p^{-\d}H)^n.
$$
\bigskip
Let us note that in Theorem A the exponent $\frac{n}{2n+2}$ tend to $1/2$ as $n\to\infty$. Theorem A was strengthened by Chang \cite{Ch}.

\bigskip

\textbf{Theorem B} (\cite{Ch}). \textit{Let $\e>0$ and a parallelepiped $B$ obeys the condition $\prod_{i=1}^nH_i > p^{(\frac25+\e)n}$. Then
$$\left|\sum_{x\in B}\chi(x)\right|\ll_{n,\e}p^{-\e^2/4}|B| 
$$
in the case where $n$ is odd and in the case where $n$ is even and $\chi|_{\F_{p^{n/2}}}$ is nontrivial character, and   
$$\left|\sum_{x\in B}\chi(x)\right|\leq \max_{\xi}|B\cap\xi\F_{p^{n/2}}|+O_{n,\e}(p^{-\e^2/4}|B|) 
$$
otherwise.}
\bigskip

Let us note that on the condition $|B|=\prod_{i=1}^nH_i>p^{(2/5+\e)n}$ it is generally impossible to obtain nontrivial estimates for sums $\sum_{x\in B}\chi(x)$ even if $\chi$ is nontrivial; indeed, one has to take into account the situation where $B$ is the subfield $\F_{p^{n/2}}$ and $\chi$ is the nontrivial character of $\F_{p^n}$ which is identical on $\F_{p^{n/2}}$. That is why one has to consider different cases which are described in Theorem B.

Further, Chang \cite{Ch2} obtained nontrivial estimates for character sums for the case $n=2$, $H_1,H_2>p^{1/4+\e}$. Konyagin \cite{Kon} generalised this result for arbitrary finite fields.

\textbf{Theorem C} (\cite{Kon}). \textit{Let $\e>0$ and $H_i > p^{1/4+\e}$ for all $1\leq i\leq n$. Then}
$$\left|\sum_{x\in B}\chi(x)\right|\ll_{n,\e}p^{-\e^2/2}|B|. 
$$

\bigskip

The aim of the present paper is to prove the following result for the cases $n=2$ and $n=3$.

\bigskip

\textbf{Theorem.} 
\textit{Let $n\in\{2,3\}$, $\chi$ be a nontrivial multiplicative character of $\F_{p^n}$ and $|B|\geq p^{n(1/4+\e)}$, and let us assume that $H_1\leq\ldots\leq H_n$. Then
$$\left|\sum_{x\in B}\chi(x)\right|\ll_{\e} |B|p^{-\e^2/12},
$$ 
if $\chi|_{\F_p}$ is not identical, and
$$\left|\sum_{x\in B}\chi(x)\right|\ll_{\e} |B|p^{-\e^2/12}+|B\cap \o_n\F_p|
$$
otherwise.} 

\bigskip

Since $\{\o_1,\ldots,\o_n\}$ is a basis, we thus have 
$$ |B\cap \o_n\F_p|=\begin{cases} H_n, \mbox{ if $0\in \cap_{i=1}^{n-1} [N_i+1,N_i+H_i]$},
\\0, \mbox{ otherwise,}
\end{cases}
$$
and in the second case of the Theorem in fact the estimate $\sum_{x\in B}\chi(x)\ll_{\e} |B|p^{-\e^2/12}+H_n$ holds. Besides, similarly to the remark for Theorem B, on the condition ($|B|\geq p^{n(1/4+\e)}$) it is generally impossible to obtain nontrivial results, since one has to keep in mind the case where $B=\F_p$ and $\chi$ is the nontrivial character which is identical on $\F_p$. Let us stress that on the condition of theorem C such a situation is impossible because of the restriction $H_i>p^{1/4+\e}$, $1\leq i\leq n$.

The key ingredient in the proofs of Theorems B and C and the Theorem of the present paper is a bound for the quantity
$$E(B)=\#\{(x,y,w,t)\in B^4 : \, xy=wt\},
$$
which is called the multiplicative energy of the set $B$. Using tools from additive combinatorics, Chang proved that $E(B)\ll_n |B|^{11/4}\log p$ for paralle-lepipeds such that $H_i<\frac12(\sqrt{p}-1)$ (see \cite{Ch}, Proposition 1 ), whereas Konyagin, using geometric number theory, established the bound $E(B)\ll_{n}|B|^2\log p$ for parallelepipeds with $H_1=\ldots=H_n\leq\sqrt{p}$ (see \cite{Kon}, Lemma 1). We generalize Lemma 1 from \cite{Kon} for the cases $n=2$, $n=3$ and distinct edges and prove the following.

\bigskip

\textbf{The Key Lemma.} \textit{Let $n\in\{2,3\}$ and suppose that
	$H_1\leq \ldots\leq H_n<\sqrt{p/2}$. Then we have}
$$E(B)\ll |B|^2\log^3 p.
$$

\bigskip

In the proof of the Theorem we closely follow \cite{Ch}: firstly, we prove the desired bound in the case where all the edges are less than $\sqrt{p/2}$ (this argument is now standard and was used in \cite{Ch}, \cite{Kon}, and had been elaborated by Karatsuba in his work \cite{Kar1}); it also immediately implies the statement for the case where all edges are less than $p^{1/2+\e/2}$. After that we prove the Theorem in the case $H_3>p^{1/2+\e/2}$. In fact, one can see from the proof that in the last case one can write a slightly better bound for the character sum, namely, $\sum_{x\in B}\chi(x)\ll_{\e} |B|p^{-\e/3}+|B\cap \o_3\F_p|$.

We prove the Key Lemma and the Theorem in the technically more difficult case $n=3$ (the case $n=2$ is absolutely similar). We prove the Key Lemma in Section 2 and the Theorem in Section 3.

\bigskip

The author would like to thank Nicholas Katz for providing an extension of his result (see Theorem E below), which is crucial for the proof of the Theorem in the case $H_3>p^{1/2+\e/2}$.

\section{Proof of the Key Lemma}

Set
$$Z'=\frac{B\setminus\{0\}}{B\setminus\{0\}}=\{z\in \F_{p^3}: \, \exists x,y\in B\setminus\{0\},\,\, xz=y \}.
$$
If $x^1,x^2,x^3,x^4\in B$, $x^1x^2=x^3x^4$ and $(x^1,x^4)\neq(0,0)$, $(x^2,x^3)\neq(0,0)$, then for some $z\in Z'$ we have $x^1z=x^3$, $x^4z=x^2$. Thus  
\begin{equation}\label{E(B)} 
E(B)\leq 2|B|^2+\sum_{z\in Z'}f^2(z),
\end{equation}
where $f(z)$ is the number of solutions to the equation $xz=y$ where $x,y\in B$. Define
$$B_0=\left\{\sum_{i=1}^3x_i\o_i \,:\, -H_i\leq x_i\leq H_i, \,\, 1\leq i\leq3\right\},
$$
$$Z=\frac{B_0\setminus\{0\}}{B_0\setminus\{0\}}, \,\,\,\, f_0(z)=\#\{(x,y)\in B_0^2 : xz=y \}.
$$
Note that if $(x_1,y_1),\ldots,(x_k,y_k)\in B^2$ are distinct solutions to the equation $xz=y$, then $(0,0),(x_2-x_1,y_2-y_1),\ldots,(x_k-x_1,y_k-y_1)$ are distinct solutions to the same equation in $B_0^2$. Thus $f(z)\leq f_0(z)$; besides, $f_0(z)=1$ for $z\in \F_{p^3}^*\setminus Z$. Therefore,
$$\sum_{z\in Z'}f^2(z)\leq \sum_{z\in Z}f_0^2(z)+|Z'\setminus Z|.
$$
Further, $|Z'|\leq |B|^2$. Recalling (\ref{E(B)}), we see that 
$$ E\leq 3|B|^2+\sum_{z\in Z} f_0^2(z),
$$
and it suffices to estimate the sum
$$S=\sum_{z\in Z}f_0^2(z).
$$
We can rewrite $S$ as
\begin{equation*} 
S=S_1+S_2,
\end{equation*}
where
\begin{equation}\label{S1}
 S_1=\sum_{z\in Z\setminus\F_p}f^2_0(z). 
\end{equation}
and
\begin{equation}\label{S2}
S_2=\sum_{z\in\F_p^*} f_0^2(z)
\end{equation}
The claim now follows from the following two lemmas.

\begin{lem}
	We have 
$$S_1\ll |B|^2\log p.
$$
\end{lem}

\begin{lem}
	We have
	$$S_2\ll |B|^2\log^3 p.
	$$
\end{lem}

\subsection{Proof of Lemma 1}

For a fixed $z\in Z$ define the lattice $\Gamma_z\subset\Z^6$:
$$\Gamma_z =\{(x_1,x_2,x_3,y_1,y_2,y_3)\in\Z^6 :  z\sum_{i=1}^3 x_i\o_i =\sum_{i=1}^3 y_i\o_i \}.
$$
For fixed $x_1,x_2,x_3\in\Z$ the condition $(x_1,x_2,x_3,y_1,y_2,y_3)\in\Gamma_z$ defines each of numbers $y_1,y_2,y_3$ modulo $p$. Thus,
\begin{multline*}
|\{(x_1,x_2,x_3,y_1,y_2,y_3)\in \Gamma_z: |x_i|,|y_i|\leq M, 1\leq i\leq 3\}|=\\
=\frac{(2M)^6}{p^3}(1+o(1)), \quad M\to\infty.
\end{multline*}
Hence
$$\mbox{\textsf{mes\,\,}}(\R^6/\Gamma_z)=p^3.
$$
Define the set 
$$D=\{(x_1,x_2,x_3,y_1,y_2,y_3)\in\R^6 : |x_i|,|y_i|\leq H_i, \,\, 1\leq i\leq3\};
$$
then we have $f_0(z)=|\Gamma_z\cap D|$. Let us recall that $i$-th successive minima 
$$\l_i=\l_i(z)=\l_i(D,\Gamma_z)
$$ 
of the set $D$ with respect to $\Gamma_z$ is defined as the least $\l>0$ such that the set $\l D$ contains $i$ linearly independent vectors of $\Gamma_z$. Obviously, $\l_1(z)\leq\ldots\leq \l_6(z)$ and $\l_1(z)\leq1$ of and only if $z\in Z$. Further, from Minkowski's second theorem (see, for instance, \cite{Tao-Vu}, Theorem 3.30) we have
\begin{equation}\label{mink} 
\prod_{i=1}^6\l_i \gg \frac{\textsf{mes}(\R^6/\Gamma_z)}{\textsf{mes}D}\gg p^3|B|^{-2}.
\end{equation}
It is well-known (see \cite{BHW}, Proposition 2.1, or Exercise 3.5.6, \cite{Tao-Vu}), that the number $f_0(z)$ of points of $\Gamma_z$ in the set $D$ obeys the inequality
\begin{equation}
f_0(z)\ll \prod_{i=1}^{6}\max\{1,\l_i^{-1}\}. \label{f_0}
\end{equation}

Now we are going to obtain lower bounds for $\l_1(z)$, $\l_2(z)$, $\l_3(z)$, where $z\in Z\setminus\F_p$.

Firstly, since $z\in Z$, then $\l_1(z)\leq1$. Besides, $H_2^{-1}\leq \l_1(z)$ (otherwise there exists a non-zero vector $(0,0,u_3,0,0,u_6)\in \Gamma_z$ such that $|u_3|,\,|u_6|<H_3H_2^{-1}$ and $zu_3\o_3=u_6\o_3$, which contradicts our assumption that $z\notin\F_p$). 

Further, we prove that $\l_2(z)\geq H_1^{-1}$. To show this, assume for contradiction that $\l_2(z)<H_1^{-1}$. Then we can find two linearly independent over $\Z$ vectors $\u=(0,u_2,u_3,0,u_5,u_6),\,\v=(0,v_2,v_3,0,v_5,v_6)\in\Gamma_z$ such that
$|u_2|,|u_5|,|v_2|,|v_5|<H_2H_1^{-1}<\sqrt{p/2}$, $|u_3|,|u_6|,|v_3|,|v_6|<H_3H_1^{-1}<\sqrt{p/2}$, and

\begin{equation}\label{system} 
\begin{cases}
(u_2\o_2+u_3\o_3)z=u_5\o_2+u_6\o_3 \,, \\
(v_2\o_2+v_3\o_3)z=v_5\o_2+v_6\o_3.
\end{cases} 
\end{equation}
	
Suppose that the vectors $(u_2,u_3)$ and $(v_2,v_3)$ are linearly independent over $\F_p$. It means that the map $x\mapsto xz$ is a bijection from the subspace $\textsf{Lin}\{\o_2,\o_3\}$ to itself. Let 
$$ z\o_1=a_1\o_1+a_2\o_2+a_3\o_3;
$$
we claim that $z=a_1$. Indeed, otherwise the map $x\mapsto x(z-a_1)$ is also a bijection from $\textsf{Lin}\{\o_2,\o_3\}$ to itself, and we have $\o_1\in\textsf{Lin}\{\o_2,\o_3\}$, which is false. Thus $z=a_1$; but that contradicts to the assumption that $z\notin\F_p$. 

Therefore the vectors $(u_2,u_3)$ and $(v_2,v_3)$ have to be linearly dependent over $\F_p$. Then the determinant of the matrix 
$\begin{pmatrix}
u_2 & u_3 \\
v_2 & v_3
\end{pmatrix}
$
equals to zero modulo $p$. But all its elements are integers bounded in magnitude by $\sqrt{p/2}$; thus the absolute value of this determinant is less than $p$, and it has to be equal to zero in $\Z$. Therefore the vectors $(u_2,u_3)$ and $(v_2,v_3)$ are linearly dependent over $\Z$. 

The vector $\v=(0,v_2,v_3,0,v_4,v_5,0)$ is non-zero; suppose that \\ $(v_2,v_3)\neq(0,0)$ and let $v_2\neq 0$ (the case $v_3\neq0$ can be easily treated in a similar way). Multiplying the second equation of (\ref{system}) by $u_2/v_2$ and subtracting it from the first one, we get
$$(u_5u_2/v_2-v_5)\o_2+(u_6u_2/v_2-v_6)\o_3=0.
$$
Since $\{\o_1,\o_2,\o_3\}$ is a basis and $|u_5u_2-v_5v_2|<p$,  $|u_6u_2-v_6v_2|<p$, then $u_5u_2-v_5v_2=u_6u_2-v_6v_2=0$, hence $\u=\frac{u_2}{v_2}\v$. But this contradicts to the fact that the vectors $\u$ and $\v$ are linearly independent over $\R$. 

Finally, if $v_2=v_3=0$, then $(v_5,v_6)\neq(0,0)$ and the same arguments are valid with $1/z$ instead of $z$; one can prove in a similar manner that the vectors $(v_5, v_6)$ and $(u_5,u_6)$ are linearly dependent and get the contradiction with the choice of $\u$ and $\v$.

\bigskip%\bigskip\bigskip\bigskip\bigskip

Thus, for $z\in Z\setminus\F_p$ we have $1\geq\l_1(z)\geq H_2^{-1}$ and $\l_3(z)\geq \l_2(z)\geq H_1^{-1}$. Define
$$Z_j=\{z\in Z\setminus\F_p: 2^{j-1}\leq H_2\l_1 <2^j \}, \,\,\, 1\leq j\leq J:=\log_2H_2+1.
$$
Note that the vector $\u\in \l_1(z)D\cap\Gamma_z$ corresponding to an element $z\in Z_j$ defines $z$. Indeed, let $\u=(u_1,u_2,u_3,u_4,u_5,u_6)$ and define the elements $x,y\in \F_{p^3}$ as follows: $x=u_1\o_1+u_2\o_2+u_3\o_3$, $y=u_4\o_1+u_5\o_2+u_6\o_3$. Then we have $z=xy^{-1}$. Therefore, $|Z_j|$ is at most the number of integers points in the box $2^jH_2^{-1}D$. Setting $j_1=\log_2(H_2/H_1)$, we see that
\begin{multline}\label{Z_j}
|Z_j|\leq |2^jH_2^{-1}D\cap\Z^6|\ll\prod_{i=1}^3\max\{1,H_i2^jH_2^{-1}\}^2\leq\\
\begin{cases}
2^{4j}H_3^2H_2^{-2},&\text{if $1\leq j<j_1$;}\\
2^{6j}|B|^2H_2^{-6},&\text{if $j_1\leq j\leq J.$}
\end{cases} 
\end{multline}

Further, set $s=s(z)=\max\{j: \l_j\leq1\}$ and 
$$Z^s=\{z\in Z\setminus\F_p :\,s(z)=s\}.$$
Recalling the definition (\ref{S1}) of the sum $S_1$, we have
$$S_1\leq \sum_{s=1}^6\sum_{z\in Z^s}f_0^2(z).
$$

\smallskip

Further we treat the cases of different $s$ in a bit routine way. 

\smallskip

For $s\leq 3$ we set $Z_j^s=Z^s\cap Z_j$. Then
$$\sum_{z\in Z^s}f_0^2(z)\leq\sum_{j} \sum_{z\in Z_j^s}f_0^2(z).
$$
We will often use the trivial bound $|Z^s_j|\leq|Z_j|$. 

Let $s=1$. By (\ref{f_0}) we have $f_0(z)\ll \l_1^{-1}$. Using (\ref{Z_j}) and the fact that for $z\in Z_j$ the bound $\l_1^{-1}(z)\ll 2^{-j}H_2$ holds, we obtain
\begin{multline} \label{s=1}
\sum_j \sum_{z\in Z_j^1}f_0^2(z)\ll \sum_{j}\sum_{z\in Z_j^1}\l_1^{-2}\ll \sum_{j}|Z_j^1|2^{-2j}H_2^2\ll \\
\sum_{1\leq j < j_1}2^{4j}H_3^2H_2^{-2}2^{-2j}H_2^2
+\sum_{j_1\leq j \leq J}|B|^22^{6j}H_2^{-6}2^{-2j}H_2^2\ll
\\
H_3^2\sum_{j\leq j_1}2^{2j}+|B|^2H_2^{-4}\sum_{j\leq J}2^{4j}\ll |B|^2. 
\end{multline}

\smallskip

Let $s=2$; by (\ref{f_0}) we have $f_0(z)\leq \l_1^{-1}\l_2^{-1}$. Let $z\in Z_j$; in the case $j<j_1$ we use the bounds $\l_1^{-1}\ll 2^{-j}H_2$ and $\l_2^{-1}\leq H_1$, 
and in the case  $j\geq j_1$ --- the bound $\l_2^{-1}\leq \l_1^{-1}\ll 2^{-j}H_2$. Also recalling (\ref{Z_j}), we see that
\begin{multline} \label{s=2}
\sum_j \sum_{z\in Z_j^2}f_0^2(z)\ll \sum_{j}\sum_{z\in Z_j^2}\l_1^{-2}\l_2^{-2}\ll \\ \sum_{1\leq j< j_1}|Z_j^2|2^{-2j}H_2^2H_1^2  + \sum_{j_1\leq j\leq J}|Z_j^2|2^{-4j}H_2^4  \ll \\
\sum_{j \leq j_1}2^{4j}H_3^2H_2^{-2}2^{-2j}H_2^2H_1^2+\sum_{j \leq J}|B|^22^{6j}H_2^{-6}2^{-4j}H_2^4 \ll  |B|^2.
\end{multline} 

\smallskip

Let $s=3$; by (\ref{f_0}) we get $f_0(z)\ll \l_1^{-1}\l_2^{-1}\l_3^{-1}$. Let $z\in Z_j$; in the case $j<j_1$ we use the bounds $\l_1^{-1}\ll 2^{-j}H_2$ and $\l_3^{-1}\leq \l_2^{-1}\leq H_1$, and in the case $j\geq j_1$ --- the bound $\l_3^{-1}\leq \l_2^{-1}\leq \l_1^{-1}\ll 2^{-j}H_2$. Again using (\ref{Z_j}), we have
\begin{multline} \label{s=3}
\sum_j \sum_{z\in Z_j^3}f_0^2(z)\ll \sum_{j}\sum_{z\in Z_j^3}\l_1^{-2}\l_2^{-2}\l_3^{-2}\ll \\
\sum_{1\leq j < j_1}2^{4j}H_3^2H_2^{-2}2^{-2j}H_2^2H_1^4 +\sum_{j_1\leq j \leq J}|B|^22^{6j}H_2^{-6}2^{-6j}H_2^6 \ll |B|^2\log p.
\end{multline} 

\smallskip

Among the cases $s>3$ we first consider $s=6$. Taking into account (\ref{mink}) and (\ref{f_0}), we obtain
\begin{equation}\label{s=6}
\sum_{z\in Z^6}f_0^2(z) \ll \sum_{z\in Z} |B|^4p^{-6} \leq |B|^2|B|^4p^{-6}\leq|B|^2p^{4\cdot 1.5-6} = |B|^2
\end{equation}
(here we use the fact that $|B|\leq p^{1.5}$, which holds due to $H_1\leq H_2\leq H_3 \leq \sqrt{p/2}$).
 
\smallskip
 
Finally, we treat the cases $s=4$ and $s=5$. Define the polar lattice $\Gamma_z^*$ as follows:
$$\Gamma_z^*=\left\{(u_1,\ldots,u_6)\in\R^6: \sum_{i=1}^3u_ix_i+\sum_{i=1}^3u_{i+3}y_i \in\Z \quad\forall (x_1,\ldots,y_3)\in \Gamma_z \right\}.
$$
Note that $\Gamma_z\supseteq p\Z^6$ implies $\Gamma_z^*\subseteq p^{-1}\Z^6$. Define the polar set 
$$D^*=\{(u_1,\ldots,v_3)\in\R^6: \sum_{i=1}^3|u_ix_i|+\sum_{i=1}^3|v_iy_i|\leq1 \,\,\text{for all}\,\,(x_1,\ldots,y_3)\in D \}
$$
Clearly
$$D^*=\{(u_1,\ldots,v_3)\in\R^6: \sum_{i=1}^3(|u_i|+|v_i|)H_i \leq1 \}.
$$
Let $\l_1^*=\l_1^*(z)$ be the first successive minima of the set $D^*$ with respect to $\Gamma_z^*$. By \cite{Ban}, Proposition 3.6, we have
\begin{equation}\label{*}
\l_1^*\l_6\ll 1.
\end{equation}
Thus, taking into account (\ref{mink}) and (\ref{f_0}), in the case $s=5$ we have
\begin{equation}\label{s=5 bound}
f_0(z)\ll \prod_{i=1}^5\l_i^{-1}(z)= \l_6(z)\prod_{i=1}^6\l_i^{-1}(z)\ll \l_6|B|^2p^{-3} \ll (\l_1^*)^{-1}|B|^2p^{-3},
\end{equation}
and in the case $s=4$ 
\begin{equation}\label{s=4 bound}
f_0(z)\ll \prod_{i=1}^4\l_i^{-1}(z)\leq \l_6^2\prod_{i=1}^6\l_i^{-1}(z)\ll \l_6^2|B|^2p^{-3} \ll (\l_1^*)^{-2}|B|^2p^{-3}.
\end{equation}

The contribution to the sum $\sum_{z\in Z^5}f_0^2(z)$ (or $\sum_{z\in Z^4}f_0^2(z)$)  from those $z\in Z^5$ (respectively $z\in Z^4$) for which $\l_1^*(z)\geq1$ can be estimated similarly to the case $s=6$ (see (\ref{s=6})). Thus we can assume $\l_1^*(z)\leq 1$. Then we have $\l_1^*\geq H_1p^{-1}$ (since if $\l<H_1p^{-1}$, then due to $\Gamma_z^*\subseteq p^{-1}\Z^6$ we see that $\l D^*\cap \Gamma_z^*=\{0\}$). Set
$$Z_j'=\{z\in Z : 2^{j-1}\leq \frac{p\l_1^*(z)}{H_1}<2^j\}, \,\, j=1,...,\log_2(p/H_1)+1\}.
$$
We claim that the vector $\u\in \l_1^*(z)D^*\cap\Gamma_z^*$ corresponding to an element $z\in Z_j'$ defines this element $z$. Suppose for contradiction that there is a non-zero vector $\u=(u_1/p,...,v_3/p)\in \Gamma_{z'}^*\cap\Gamma_{z''}^*$, where $z'\neq z''$ and $u_i,v_i\in\Z$; we also have $\sum_{i=1}^3|u_i|+\sum_{i=1}^3|v_i|<2^j$. Take an arbitrary element $x=\sum_{i=1}^3x_i\o_i\in\F_{p^3}$ and set $y'=xz'=\sum_{i=1}^3y'_i\o_i$, $y''=xz''=\sum_{i=1}^3y''_i\o_i$. Then 
$$(x_1,x_2,x_3,y_1',y_2',y_3')\in\Gamma_{z'}, \,\, (x_1,x_2,x_3,y_1'',y_2'',y_3'')\in\Gamma_{z''},
$$
and by the definition of the polar set
$$\sum_{i=1}^3x_iu_i/p+\sum_{i=1}^3y_i'v_i/p\in\Z, \,\, \sum_{i=1}^3x_iu_i/p+\sum_{i=1}^3y_i''v_i/p\in\Z.
$$
But then
$$\sum_{i=1}^3(y_i'-y_i'')v_i\equiv 0 {\pmod p}.
$$
Note the numbers $y_i'-y_i''$ can be arbitrary (they are the coefficients of the element $y'-y''$ which is equal to $x(z'-z'')$ and, since $z'-z''\neq0$, can be equal to a given element provided we take the appropriate $x$). Thus $v_i\equiv0\pmod{p}$, and since $|v_i|<2^j\leq p$, then $v_i=0$. So we see that
$$\sum_{i=1}^3x_iu_i=0
$$
for all $(x_1,x_2,x_3)\in\Z^3$, and hence $u_i=0$. But this contradicts to the fact that the vector $(u_1,\ldots,v_3)$ is non-zero. Therefore, the vector $\u\in \l_1^*(z)D^*\cap\Gamma_z^*$ corresponding to an element $z\in Z_j'$ indeed defines $z$.

The vector $(u_1,\ldots,v_3)\in \frac{2^jH_1}{p}D^*\cap\Gamma_z^*$ obeys the inequality $\sum_{i=1}^3(|u_i|+|v_i|)H_i\leq 2^jH_1$; hence, $|u_i|,|v_i|\leq 2^jH_1H_i^{-1}$. Thus we see that $|Z_j'|\leq \prod_{i=1}^3\max(1,2^jH_1H_i^{-1})^2$. Setting $j_2=\log_2(H_3/H_1)$ and $j_3=\log_2(p/H_1)+1$, we have  
\begin{equation}
\label{Z_j'}
|Z'_j|\leq\begin{cases}
2^{2j}, \mbox{ if $1\leq j < j_1;$}\\
2^{4j}H_1^2H_2^{-2}, \mbox{ if $j_1\leq j < j_2$;}\\
2^{6j}H_1^6|B|^{-2}, \mbox{ if $ j_2 \leq j\leq j_3$.}\\
\end{cases}
\end{equation}
For $s=4$ and $s=5$ define
$$Z^s_j=Z^s\cap Z_j';
$$
below we will use the trivial bound $|Z^s_j|\leq |Z'_j|$ and apply (\ref{Z_j'}). Recalling the bound (\ref{s=5 bound}) and taking into account that $\l_1^*(z)\asymp 2^jH_1/p$ for $z\in Z'_j$, we obtain
\begin{multline}\label{s=5}
\sum_{z\in Z^5}f_0^2(z)\leq |B|^4p^{-6}\sum_{j}\sum_{z\in Z^5_j}(\l_1^*(z))^{-2} \leq \\
|B|^4p^{-6}\sum_{1\leq j<j_1}2^{2j-2j}H_1^{-2}p^2+ 
|B|^4p^{-6}\sum_{j_1\leq j< j_2}H_1^2H_2^{-2}2^{4j-2j}H_1^{-2}p^2+\\
|B|^4p^{-6}\sum_{j_2\leq j\leq j_3}H_1^6|B|^{-2}2^{6j-2j}H_1^{-2}p^2 \leq \\
|B|^4p^{-4}H_1^{-2}\log p +|B|^4p^{-4}H_2^{-2}\sum_{j\leq j_2}2^{2j} + |B|^2p^{-4}H_1^4\sum_{j\leq j_3}2^{4j} \ll \\
|B|^2(1+|B|^2p^{-4}H_3^4|B|^{-2}+|B|^2p^{-4}H_1^{-2}\log p)\ll |B|^2.
\end{multline} 

In the case $s=4$, using (\ref{s=4 bound}), in a similar way we get
\begin{multline}\label{s=4}
\sum_{z\in Z^5}f_0^2(z)\leq |B|^4p^{-6}\sum_{j}\sum_{z\in Z^5_j}(\l_1^*(z))^{-4} \leq \\
|B|^4p^{-6}\sum_{1\leq j<j_1}2^{2j-4j}H_1^{-4}p^4+
|B|^4p^{-6}\sum_{j_1\leq j< j_2}H_1^2H_2^{-2}2^{4j-4j}H_1^{-4}p^4+\\
|B|^4p^{-6}\sum_{j_2\leq j\leq j_3}H_1^6|B|^{-2}2^{6j-4j}H_1^{-4}p^4 \leq\\
|B|^4p^{-2}H_1^{-4}+|B|^4p^{-2}H_1^{-2}H_2^{-2}\sum_{j\leq j_2}1 + |B|^2p^{-2}H_1^2\sum_{j\leq j_3}2^{2j} \ll \\
|B|^2(1+|B|^2p^{-2}H_3^2(\log p)|B|^{-2}+|B|^2p^{-2}H_1^{-4})\ll |B|^2.
\end{multline} 
\bigskip
Putting the bounds (\ref{s=1})-(\ref{s=6}) and (\ref{s=5})-(\ref{s=4}) together, we see that
\begin{equation*}
S_1\ll |B|^2\log p,
\end{equation*} 
as desired.

\subsection{Proof of Lemma 2}

Fix $z\in\F_p$. Let $x=\sum_{i=1}^3x_i\o_i$ and $y=\sum_{i=1}^3y_i\o_i$; then the equality $xz=y$ is equivalent to the equalities $zx_i\equiv y_i\pmod{p}$, $1\leq i\leq3$. Hence
$$f_0(z)=f_1(z)f_2(z)f_3(z),
$$
where
$$f_i(z)=\#\{(x_i,y_i)\in[-H_i,H_i]^2 : x_iz\equiv y_i \pmod{p} \}.
$$
Recalling the definition (\ref{S2}) of the sum $S_2$, we see that
\begin{equation} \label{S_2 bound}
S_2=\sum_{z\in \F_p^*}f_0^2(z)=\sum_{z\in\F_p}f_1^2(z)f_2^2(z)f_3^2(z)\leq\prod_{i=1}^3\left(\sum_{z\in \F_p^*}f_i^2(z)\right).
\end{equation} 
The sums $\sum_{z\in \F_p^*}f_i^2(z)$ can be estimated as the sum $S_1$ in the previous subsection. We go over the details quickly. Fix $i\in\{1,2,3\}$ and denote for the brevity $H=H_i$,
$$D=[-H,H]^2,  \qquad Z=\frac{[-H,H]\setminus{0}}{[-H,H]\setminus{0}}\,,
$$
$$\Gamma_z=\{(x,y)\in\Z^2 : \,\,xz\equiv y \pmod{p}\};
$$
let $\l_l=\l_l(z)$ be the $l$-th successive minima of $D$ with respect to $\Gamma_z$, $l=1,2$. Then for all $z\in\F_p^*$ we have
$$\textsf{mes}(\R^2/\Gamma_z)=p,
$$
and Minkowski's second theorem gives us 
\begin{equation}\label{mink2}
 \l_1\l_2\gg pH^{-2}.
\end{equation} 
In our notation we have $f_i(z)=|D\cap\Gamma_z|$. By Proposition 2.1 from \cite{BHW} we see that
$$
f_i(z)\ll \prod_{l=1}^2\max\{1,\l_l^{-1}(z)\}.
$$  
Clearly, $H^{-1}\leq\l_1\leq 1$ for $z\in Z$. Define the set 
$$Z_j=\{z\in Z: 2^{j-1}\leq H\l_1(z) < 2^j \}, \,\, j=1,\ldots,[\log_2H]+1,
$$
and let $s(z)=\max\{l: \l_l(z)\leq 1\}$ and $Z^s=\{z\in Z : \,\,s(z)=s\}$. The vector $(u_1,u_2)\in \l_1(z)D\cap\Gamma_z$ corresponding to an element $z\in Z_j$ defines $z$. Thus
$$|Z_j|\ll \left|\frac{2^j}{H}D\cap \Z^2\right| \ll 2^{2j}
$$
and
\begin{multline}\label{S_2, s=1}
\sum_{z\in Z^1}f_i^2(z) \ll \sum_{j}\sum_{z\in Z^1\cap Z_j} \l_1^{-2}(z) \ll \sum_{j=1}^{[\log_2H]+1}2^{2j}H^22^{-2j}\ll H^2\log p. 
\end{multline}
Finally, using (\ref{mink2}) and the fact that $H\leq \sqrt{p}$, we find
\begin{equation}\label{S_2,s=2}
\sum_{z\in Z^2}f_i^2(z)\ll |Z^2|p^{-2}H^4 \leq H^6p^{-2}\leq H^2. 
\end{equation}
Putting (\ref{S_2, s=1}) and (\ref{S_2,s=2}) together, we obtain
$$\sum_{z\in\F^*_p}f_i^2(z)\ll H^2\log p.
$$
Recalling (\ref{S_2 bound}), we get
$$
S_2 \ll \prod_{i=1}^3(H_i^2\log p)=|B|^2\log^3 p.
$$
This completes the proof of Lemma 2 and the Key Lemma.

\section{Proof of the Theorem.} 

In this section we closely follow to the paper \cite{Ch}. We would like to stress, however, that the arguments in the case $H_3<\sqrt{p/2}$ (additive shift $x\mapsto x+yz$ and double application of H\"older's inequality) are now standard and were used in works \cite{Ch}, \cite{Kon} and had been elaborated by Karatsuba in his work \cite{Kar1}. Additive shift itself was used earlier in works of Vinogradov (see \cite{Vin1}, \cite{Vin2}, \cite{Vin3}) and probably rises from ideas of van der Corput and H.Weil (see, for instance, \cite{vdC}, \cite{Weil1}, \cite{Weil2}).

\subsection{The case $H_3 < \sqrt{p/2}$.}

Dividing $B$ to smaller parallelepipeds, we may assume that $|B|\asymp p^{3(1/4+\e)}$. Let $\d=\d(\e)>0$ be chosen later. Set
$$I=[1,p^{\d}]\cap\Z 
$$
and
$$B_0=\left\{\sum_{i=1}^3x_i\o_i : \,\, x_i\in[0,p^{-2\d}H_i]\cap\Z, \,\, 1\leq i\leq3 \right\}.
$$
Note that $\#([0,p^{-2\d}H_i]\cap\Z)\asymp 1+p^{-2\d}H_i\gg p^{-2\d}H_i$, and, hence, we have
\begin{equation}\label{B_0}
|B_0|\gg p^{-6\d }|B|.
\end{equation}
Since $B_0I\subseteq \left\{\sum_{i=1}^3x_i\o_i: \,\, x_i\in[0,p^{-\d}H_i]\cap\Z, \,\, 1\leq i\leq3 \right\}$, for all $y\in B_0$, $z\in I$ we have 
$$\left|\sum_{x\in B}\chi(x)-\sum_{x\in B}\chi(x+yz)\right|\leq |B\setminus(B+yz)|+|(B+yz)\setminus B| \leq 6p^{-\d}|B|. 
$$
Thus
\begin{equation}\label{b1} 
\sum_{x\in B}\chi(x)=\frac{1}{|B_0||I|}\sum_{x\in B, y\in B_0, z\in I}\chi(x+yz)+O(p^{-\d}|B|). 
\end{equation}
Further,
\begin{multline*}
\left|\sum_{x\in B,y\in B_0,z\in I}\chi(x+yz)\right|\leq \sum_{x\in B,y\in B_0}\left|\sum_{z\in I}\chi(x+yz)\right| \leq \\
\sum_{x\in B,y\in B_0\setminus\{0\}}\left|\sum_{z\in I}\chi(xy^{-1}+z)\right|+ |B||I|=\\
 \sum_{u\in \F_{p^3}}\tau(u)\left|\sum_{z\in I}\chi(u+z)\right|+ |B||I|,
\end{multline*}
where
$$\tau(u)=\#\{(x,y)\in B\times (B_0\setminus\{0\}) : xy^{-1}=u\}.
$$
Let $r$ be a positive integer to be chosen later. Using H\"older's inequality twice, we obtain
\begin{multline}\label{b2}
\left|\sum_{x\in B,y\in B_0,z\in I}\chi(x+yz)\right|\leq \\ \left(\sum_{u\in\F_{p^3}}\tau(u)\right)^{1-1/r}\left(\sum_{u\in\F_{p^3}}\tau(u)\left|\sum_{z\in I}\chi(u+z)\right|^r\right)^{1/r}+|B||I| \leq \\
 \left(\sum_{u\in\F_{p^3}}\tau(u)\right)^{1-1/r}\left(\sum_{u\in\F_{p^3}}\tau^2(u)\right)^{1/(2r)}\left(\sum_{u\in \F_{p^3}}\left|\sum_{z\in I}\chi(u+z)\right|^{2r}\right)^{1/(2r)}+|B||I|. 
\end{multline}
Now we have to estimate three sums which have appeared in the last line of (\ref{b2}).
Firstly, 
\begin{equation}\label{o}
 \sum_{u\in\F_{p^3}}\tau(u) =|B|(|B_0|-1)\leq |B||B_0|.
\end{equation}
Further, $\tau(0)\leq |B_0|$ and hence
$$\tau(0)^2\leq |B_0|^2\leq |B||B_0|.
$$
Using the Cauchy-Schwarz inequality and the Key Lemma, we see that
\begin{multline*}
\sum_{u\in\F_{p^3}^*}\tau^2(u) =\#\{(x_1,x_2,y_1,y_2)\in B\times B \times B_0\times B_0\,:\,x_1y_2=x_2y_1\neq 0 \}=\\
\sum_{\nu\in\F_{p^3}^*}\#\{(x_1,x_2)\in B^2:\frac{x_1}{x_2}=\nu \}\#\{(y_1,y_2)\in B_0^2:\frac{y_1}{y_2}=\nu \}\leq \\
E(B)^{1/2}E(B_0)^{1/2}\ll 
|B||B_0|\log^3p.
\end{multline*}
Putting together the last two inequalities, for the second sum we get the bound
\begin{equation}\label{o^2}
\sum_{u\in\F_{p^3}}\tau^2(u) \ll |B||B_0|\log^3 p.
\end{equation}

In order to estimate the third sum we will use the following theorem.

\bigskip

\textbf{Theorem D} (\cite{Sch}, Theorem 2C', p.43). \textit{Let $\chi$ be a multiplicative character of $\F_{p^n}$ of order $d>1$. Assume that a polynom $f\in\F_{p^n}[x]$ has $m$ distinct roots  and is not $d$-th power. Then} 
$$\left|\sum_{x\in \F_{p^n}}\chi(f(x))\right|\leq (m-1)p^{n/2}.
$$

\bigskip

We have
\begin{multline*}
\sum_{u\in\F_{p^3}}\left|\sum_{z\in I}\chi(u+z)\right|^{2r}\leq\\
\sum_{z_1,\ldots,z_{2r}\in I}\left|\sum_{u\in\F_{p^3}}\chi\left((u+z_1)\ldots(u+z_r)(u+z_{r+1})^{q-2}\ldots(u+z_{2r})^{q-2}\right)\right|.
\end{multline*}
We call a tuple $(z_1,\ldots,z_{2r})$ \textit{good} if at least one of its elements occurs exactly once, and call it \textit{bad} otherwise. By Theorem D we have the bound
$$\left|\sum_{u\in\F_{p^3}}\chi\left((u+z_1)\ldots(u+z_r)(u+z_{r+1})^{q-2}\ldots(u+z_{2r})^{q-2}\right)\right| < 2rp^{3/2}
$$
for any good tuple $(z_1,\ldots,z_{2r})$. We can estimate the number of good tuples trivially by $|I|^{2r}$ and thus see that the contribution from them is at most $2rp^{3/2}|I|^{2r}$. Further, in any bad tuple every element occurs at least twice, and hence it contains at most $r$ distinct element. They can be chosen in at most $|I|^r$ ways, and hence the number of bad tuples does not exceed $|I|^rr^{2r}$. We can estimate the contribution from each bad tuple trivially by $p^3$, and thus see that the contribution from bad tuples is at most $p^3|I|^rr^{2r}$. Therefore,
\begin{equation*}
\sum_{u\in\F_{p^3}}\left|\sum_{z\in I}\chi(u+z)\right|^{2r}\leq 2rp^{3/2}|I|^{2r}+p^3|I|^rr^{2r},
\end{equation*}
and hence
\begin{equation}\label{o^3}
\left(\sum_{u\in\F_{p^3}}\left|\sum_{z\in I}\chi(u+z)\right|^{2r}\right)^{1/(2r)}\ll p^{3/(4r)}|I|+p^{3/(2r)}|I|^{1/2}r.
\end{equation}
Putting the bounds (\ref{o})-(\ref{o^3}) into (\ref{b2}), we get
\begin{multline*}
\frac{1}{|B_0||I|}\left|\sum_{x\in B,y\in B_0,z\in I}\chi(x+yz)\right|\ll \\
\frac{1}{|B_0||I|}\left(|B||B_0|\right)^{1-\frac{1}{r}}\left(|B||B_0|\log^3p\right)^{\frac{1}{2r}}\left(p^{3/(4r)}|I|+rp^{3/(2r)}|I|^{1/2}\right)+|B||B_0|^{-1}\\
=|B|\left(|B||B_0|\right)^{-1/(2r)}(\log p)^{3/(2r)}\left(p^{3/(4r)}+rp^{3/(2r)}|I|^{-1/2}\right)+|B||B_0|^{-1}.
\end{multline*}
Recalling the bound (\ref{B_0}) and the assumption on the quantity $|B|$ and taking into account the $|I|\gg_{\e} p^{\d}/2$ (recall that $\d$ will be depending only on $\e$), we have
\begin{multline*}
\frac{1}{|B_0||I|}\left|\sum_{x\in B,y\in B_0,z\in I}\chi(x+yz)\right|\ll_{\e} \\
|B|p^{-3/(4r)-3(\e-\d)/r}(\log p)^{3/(2r)}\left(p^{3/(4r)}+rp^{3/(2r)-\d/2}\right)+O(p^{6\d}).
\end{multline*}
Set $\d=\frac{3}{2r}$. Then
\begin{equation*}
\frac{1}{|B_0||I|}\left|\sum_{x\in B,y\in B_0,z\in I}\chi(x+yz)\right|\ll_{\e}
|B|rp^{-3(\e-\d)/r}(\log p)^{3/(2r)}+O(p^{6\d}).
\end{equation*}
Recalling (\ref{b1}) and the fact that $3/r=2\d$, we get
\begin{equation}\label{final bound1}
\left|\sum_{x\in B}\chi(x)\right| \ll_{\e} |B|rp^{-2\d(\e-\d)}(\log p)^{\d}+p^{6\d}+|B|p^{-\d}.
\end{equation}
We choose $r$ so that $\d=3/(2r)$ is close to $\e/2$. To be more precise, let $r$ be the nearest integer to the number $3\e^{-1}$; then
$$
\left|r-\frac{3}{\e}\right|\leq1/2 
$$
and
$$r=3\e^{-1}+0.5\theta,
$$
where $|\theta|\leq1$. Thus
$$\d=\frac{3}{2r}=\frac{3}{2(3\e^{-1}+0.5\theta)}=\frac{\e}{2+\theta\e/3}
$$
and hence $\frac13\e<\frac{6}{13}\e\leq \d \leq \frac{6}{11}\e$ (we may assume $\e<1/2$). Since $|B|\gg p^{3/4+3\e}$, then $p^{6\d}\ll |B|p^{-\d}\leq|B|p^{-\e/3}$, and we can rewrite (\ref{final bound1}) as
\begin{equation*}\label{final bound2}
\left|\sum_{x\in B}\chi(x)\right| \ll_{\e} |B|p^{-2\d(\e-\d)}(\log p)^{\d}+|B|p^{-\e/3}.
\end{equation*}
Finally,
$$2\d(\e-\d)\geq 2(6\e/13)(5\e/11)=60\e^2/143>\e^2/3.
$$
Hence
$$\left|\sum_{x\in B}\chi(x)\right| \ll_{\e} |B|p^{-\e^2/3}. 
$$
This concludes the proof of the Theorem in the case $H_3\leq\sqrt{p/2}$.

\subsection{The case $H_3\leq p^{1/2+\e/2}$.}

In this case we can divide each edge which has length greater than $\sqrt{p/2}$ into $O(p^{\e/2})$ ``almost equal'' pieces of length less than $\sqrt{p/2}$ but greater than $\sqrt{p}/2$. So $B$ can be divided into $O((p^{\e/2})^3)$ parallelepipeds $B_{\alpha}$ of volume $\gg p^{-3\e/2}p^{3(1/4+\e)}= p^{3(1/4+\e/2)}$. According to the previous case
$$\left|\sum_{x\in B_{\alpha}}\chi(x)\right|\ll_{\e} |B_{\alpha}|p^{-\e^2/12}
$$
for all $\alpha$, and thus 
$$\left|\sum_{x\in B}\chi(x)\right|\ll_{\e} |B|p^{-\e^2/12}.
$$

\subsection{The case $H_3>p^{1/2+\e/2}$.}

We need the following extension of a result of Katz \cite{K}.

\bigskip 

\textbf{Theorem E}. \textit{Let $\chi$ be a nontrivial multiplicative character of $\F_{p^n}$ and $g\in\F_{p^n}$ be a generating element, i.e., $\F_{p^n}=\F_p(g)$. Then for any interval $I\subseteq [1,p]\cap\Z$ we have}
$$\left|\sum_{t\in I}\chi(g+t)\right|\leq c(n)\sqrt{p}\log p.
$$

\bigskip 

We can rewrite the initial sum as
\begin{equation}\label{3.1}
\left|\sum_{x\in B}\chi(x)\right|=\left|\sum_{(x_1,x_2)\in I_1\times I_2}\sum_{x_3\in I_3}\chi(x_1\frac{\o_1}{\o_3}+x_2\frac{\o_2}{\o_3}+x_3)\right|,
\end{equation}
where $I_i=[N_i+1,N_i+H_i]\cap\Z$. Define the set $A$ as follows: 
$$A=\left\{(x_1,x_2)\in I_1\times I_2 \,:\,\F_p(x_1\frac{\o_1}{\o_3}+x_2\frac{\o_2}{\o_3})\neq\F_{p^3} \right\}.
$$
Since $3$ is a prime number, $\F_p$ is the only nontrivial subfield of $\F_{p^3}$, and we have
$$A=\left\{(x_1,x_2)\in I_1\times I_2 \,:\,x_1\frac{\o_1}{\o_3}+x_2\frac{\o_2}{\o_3}\in\F_p \right\}.
$$
Further, the elements $1,\frac{\o_1}{\o_3},\frac{\o_2}{\o_3}$ are linearly independent over $\F_p$, and hence $x_1\frac{\o_1}{\o_3}+x_2\frac{\o_2}{\o_3}\in\F_p$ if and only if $x_1=x_2=0$. We thus see that
$$A=\begin{cases}
\{(0,0)\}, \mbox{ if $0\in I_1\cap I_2.$}\\
\varnothing, \mbox{ otherwise.}
\end{cases}
$$
Now let us turn to equality (\ref{3.1}). If a pair $(x_1,x_2)$ does not belong to $A$, then by Theorem E and the assumption on $H_3$ we have
$$\left|\sum_{x_3\in I_3}\chi(x_1\frac{\o_1}{\o_3}+x_2\frac{\o_2}{\o_3}+x_3)\right|\ll \sqrt{p}\log p \leq H_3p^{-\e/2}\log p.
$$
Thus we can bound the number of pairs $(x_1,x_2)$ which do not belong to $A$ trivially by $|I_1||I_2|$, we obtain 
$$\left|\sum_{(x_1,x_2)\in (I_1\times I_2)\setminus D}\sum_{x_3\in I_3}\chi(x_1\o_1+x_2\o_2+x_3\o_3)\right| \ll_{\e} |B|p^{-\e/3}. 
$$
This concludes the proof of the Theorem in the case $0\notin I_1\cap I_2$. Now suppose that $0\in I_1\cap I_2$. By arguing as before we see that it suffices to estimate the sum
$$S'=\sum_{x_3\in I_3}\chi(x).
$$
If $\chi|_{\F_p}$ is not identical, then by the Polya-Vinogradov inequality and the assumption on $H_3$ we have
$$|S'|\leq \sqrt{p}\log p \ll_{\e} H_3p^{-\e/3}\ll_{\e} |B|p^{-\e/3}.
$$
This completes the proof in the case where  $0\in I_1\cap I_2$ and $\chi|_{\F_p}$ is not identical.

Finally we consider the case where $\chi|_{\F_p}$ is the trivial character. Then
$$|S'|\leq H_3,
$$
and thus we see that in the case $H_3\geq p^{1/2+\e/2}$ we always have the bound
$$\left|\sum_{x\in B}\chi(x)\right|\ll_{\e} |B|p^{-\e/3}+H_3.
$$

The claim follows.

\end{document}